\documentclass[12pt]{article}
\usepackage{amssymb}
\usepackage{mathrsfs}
\usepackage[hypertex]{hyperref}
\usepackage{amsfonts}
\usepackage{graphicx}
\usepackage{mathptmx}
\usepackage{latexsym,amsmath,amssymb,amsfonts,amsthm}
\newcommand{\R}{\mathbb{R}}

\newcommand{\bcen}{\begin{center}}
\newcommand{\ecen}{\end{center}}
\newtheorem{theorem}{Theorem}[section]

\newtheorem{corollary}[theorem]{Corollary}

\newtheorem{remark}[theorem]{Remark}

\newtheorem{defn}[theorem]{Definition}

\def\inx{\int_X}
\def\inr{\int_{\mathbb{R}^n}}
\def\ino{\int_0^{+\infty}}
\def\inm{\int_M}
\def\dn{d\nu_g}
\def\dm{d\mu(x)}
\def\de{d\mu_E(x)}
\def\n{\nabla}
\def\la{\lambda}

\def\p{\partial}
\def\R{\mathbb{R}}
\def\({\left(}
\def\){\right)}

\begin{document}
\setcounter{page}{1}
\title{The Gagliardo-Nirenberg inequality on metric measure
spaces}
\author{Feng Du$^{a,b}$, ~Jing Mao$^{c,d,\ast}$,~
Qiaoling Wang$^{b}$,~ Chuanxi Wu$^e$}

\date{}
\protect\footnotetext{\!\!\!\!\!\!\!\!\!\!\!\!{$^{\ast}$Corresponding author}
{{\emph{E-mail addresses}}: defengdu123@163.com(F. Du),
jiner120@163.com(J. Mao), wang@mat.unb.br(Q. Wang),
cxwu@hubu.edu.cn(C. Wu).}}
\maketitle ~~~\\[-15mm]

\begin{center}
{\footnotesize  $a$. School of Mathematics and Physics Science,
Jingchu University of Technology, Jingmen, 448000, China\\
$b$. Departamento de Matem\'atica, Universidade de Brasilia,
70910-900-Brasilia-DF, Brazil\\
$c$. Department of Mathematics, Harbin Institute of Technology (Weihai), Weihai, 264209, China\\
$d$. Instituto Nacional de Matem\'{a}tica Pura e Aplicada, 110
Estrada Dona Castorina, Rio de Janeiro, 22460-320, Brazil\\
$e$. School of Mathematics and Statistics, Hubei University, Wuhan,
430062, China }
\end{center}


\begin{abstract}
In this paper, we prove that if a metric measure space satisfies the
volume doubling condition and the Gagliardo-Nirenberg inequality
with the same exponent $n$ $(n\geq 2)$, then it has exactly the
$n$-dimensional volume growth. Besides, two interesting applications
have also been given. The one is that we show that if a complete
$n$-dimensional Finsler manifold of nonnegative $n$-Ricci curvature
satisfies the Gagliardo-Nirenberg inequality with the sharp
constant, then its flag curvature is identically zero. The other one
is that we give an alternative proof to Mao's main result in
\cite{m3} for smooth metric measure spaces with nonnegative weighted
Ricci curvature.
 \end{abstract}

\markright{\sl\hfill F. Du, J. Mao, Q.-L. Wang, C-X. Wu  \hfill}

\section{Introduction}
\renewcommand{\thesection}{\arabic{section}}
\renewcommand{\theequation}{\thesection.\arabic{equation}}
\setcounter{equation}{0} \label{intro} Let $M$ be an $n$-dimensional
complete non-compact Riemannian manifold, and denote by $\nabla$ the
gradient operator on $M$. Given positive numbers $p$ and $q$, denote
by $\mathcal{D}^{p,q}(M)$ the completion of the space of smooth
compactly supported functions on $M$ under the norm$\|\cdot\|_{p,q}$
defined by $\|u\|_{p,q}=\|\nabla u\|_p+\|u\|_q$ . Let $1<p<n$,
$p<q\leq \frac{p(n-1)}{n-p}$, $\delta=np-(n-p)q$,
$r=p\frac{q-1}{p-1}$, $\theta=\frac{(q-p)n}{(q-1)(np-(n-p)q)}$. For
all $u\in \mathcal{D}^{p,q}(\mathbb{R}^n)$, Del Pino-Dolbeault
\cite{dpd1,dpd2} proved that
\begin{eqnarray}\label{1.1}
\(\inr |u|^rdx\)^{\frac{1}{r}}\leq\Phi\(\inr |\nabla
u|^pdx\)^{\frac{\theta}{p}}\(\inr |u|^qdx\)^{\frac{1-\theta}{q}},
 \end{eqnarray}
where $\Phi$ is the best constant for the inequality (\ref{1.1}) and
takes the explicit form
\begin{eqnarray*}
\Phi=\(\frac{q-p}{p\sqrt{\pi}}\)^\theta\(\frac{pq}{n(q-p)}\)^{\frac{\theta}{p}}
\(\frac{\theta}{pq}\)^{\frac{1}{r}}\(\frac{\Gamma\(q\frac{p-1}{q-p}\Gamma\(\frac{n}{2}+1\)\)}
{\Gamma\(\frac{p-1}{p}\frac{\delta}{q-p}\)\Gamma\(n\frac{p-1}{p}+1\)}\)^{\frac{\theta}{n}}.
 \end{eqnarray*}
Equality holds in (\ref{1.1}) if and only if for some $\alpha\in
\mathbb{R}$, $\beta >0$, $\overline{x} \in \mathbb{R}^n$,
\begin{eqnarray*}
u(x)=\alpha\(1+\beta|x-\overline{x}|^{\frac{p}{p-1}}\)^{-\frac{p-1}{q-p}},
\qquad \forall x\in \mathbb{R}^n.
\end{eqnarray*}
The inequality (\ref{1.1}) is usually called the
\emph{Gagliardo-Nirenberg inequality}. Moreover, when
$q=p\frac{n-1}{n-p}$, then $\theta=1$, $r=\frac{np}{n-p}$, and
(\ref{1.1}) becomes the optimal Sobolev inequality, which is
separately found by Aubin \cite{au1} and Telenti \cite{t}, having
many important applications (see, for instance,
\cite{au2,au3,hlp,h1,l}). Complete manifolds with nonnegative Ricci
curvature on which some Sobolev or Caffarelli-Kohn-Nirenberg type
inequality is satisfied were studied in \cite{dcx,le,x2}.

Let $M$ be a Riemannian manifold. Let $dv$ be the Riemannian volume
element on $M$, and $C_0^\infty(M)$ be the space of smooth functions
on $M$ with compact support. Let $B(x, r)$ be the geodesic ball with
center $x \in M$ and radius $r$, and $\mathrm{Vol}[B(x, r)]$ be the
volume of $B(x, r)$, which is given by
\begin{eqnarray*}
\mathrm{Vol}[B(x, r)]=\int_{B(x, r)}dv.
\end{eqnarray*}
 In 2005, Xia
\cite{x1} studied complete non-compact Riemannian manifolds with
nonnegative Ricci curvature on which some Gagliardo-Nirenberg type
inequality is satisfied, and proved the following result.

\begin{theorem}  \label{theorem1.1}
Let $1<p<n, p<q\leq
\frac{p(n-1)}{n-p},r=p\frac{q-1}{p-1},\theta=\frac{(q-p)n}{(q-1)(np-(n-p)q)}$,
and let $C\geq\Phi$ be a constant. Assume that $M$ is an
$n$-dimensional $(n\geq2)$ complete non-compact Riemannian manifold
with non-negative Ricci curvature and assume that for any $u \in
C_0^\infty(M)$, we have
\begin{eqnarray}\label{1.2}
\(\inm |u|^rdv\)^{\frac{1}{r}}\leq C\(\inm|\nabla
u|^pdv\)^{\frac{\theta}{p}}\(\inm |u|^qdv\)^{\frac{1-\theta}{q}}.
 \end{eqnarray}
Then for any $x \in M$, we have
\begin{eqnarray}\label{1.3}
\mathrm{Vol}[B(x,r)]\geq\(C^{-1}\Phi\)^{\(\frac{\theta}{p}+\frac{1-\theta}{q}-\frac{1}{r}\)^{-1}}V_0(r),\qquad
\forall r>0,
 \end{eqnarray}
 where $V_0(r)$ is the volume of an $r$-ball in $\mathbb{R}^n$.
 \end{theorem}

Let $(X,d)$ be a metric measure space, and $\mu$ be a Borel measure
on $X$ such that $0<\mu (U)<\infty$ for any nonempty bounded open
set $U\subset X$. Let $\mathrm{Lip_0}(X)$ be the space of Lipschitz
functions with compact support on $X$,
 and define $|D u|(x)$ as follows
 \begin{eqnarray*}
 |D u|(x):=\limsup_{y\rightarrow x}\frac{|u(y)-u(x)|}{d(x,y)},
 \end{eqnarray*}
 which is the local Lipschitz constant of $u$ at $x \in X$. The
 function $x\rightarrow|D u|(x)$ is Borel measurable for
 $u\in\mathrm{Lip_0}(X)$.
 In 2013, Krist\'aly-Ohta \cite{ko} studied metric measure spaces satisfying the volume doubling
condition mentioned therein and the Caffarelli-Kohn-Nirenberg
inequality with the same exponent $n \geq 3$, and then they proved
that those spaces have exactly the $n$-dimensional volume growth.
Inspired by Xia's and Krist\'aly-Ohta's works mentioned above, here
we investigate a metric measure space satisfying a volume doubling
condition and the Gagliardo-Nirenberg inequality, and successfully
prove the following result.

 \begin{theorem} \label{theorem1.2}
 Let $p$, $q$, $r$, $\theta$, $n$ be as in Theorem \ref{theorem1.1}, $x_0\in X$, $C\geq \Phi$, and $C_0\geq 1$.  Assume that for any $u\in
\mathrm{Lip_0}(X)$, the Gagliardo-Nirenberg inequality
\begin{eqnarray}\label{a1}
\left(\inx |u(x)|^r\dm\right)^{\frac{1}{r}}\leq C\left(\inx |\n
u|^p(x)\dm\right)^{\frac{\theta}{p}}\left(\inx
|u(x)|^q\dm\right)^{\frac{1-\theta}{q}}
 \end{eqnarray}
  and the volume conditions
 \begin{eqnarray}\label{a2}
\frac{\mu (B(x,R))}{\mu (B(x,r))}\leq
C_0\left(\frac{R}{r}\right)^{n}, \quad \mathrm{for~all}~ x\in X,
~\mathrm{and}~ 0<r<R,
\end{eqnarray}
 \begin{eqnarray}\label{a3}
 \liminf_{r\rightarrow 0}\frac{\mu (B(x_0,r))}{\mu_{E} (\mathbb{B}_{n}(r))}=1
\end{eqnarray}
hold on a proper metric measure space $(X,d,\mu)$ of dimension $n$,
where $B(x, r ) := \{y \in X : d(x, y) < r \}, \mathbb{B}_n(r ) :=
\{x
 \in \mathbb{R}^n : |x| < r \}$, and $\mu_E$
is the $n$-dimensional Lebesgue measure. Then, for any $x\in X$ and
$\rho>0$, we have
\begin{eqnarray}\label{a4}
 \mu (B(x,\rho))\geq C_{0}^{-1}\(C^{-1}\Phi\)^{\(\frac{\theta}{p}+\frac{1-\theta}{q}-\frac{1}{r}\)^{-1}}\mu_E(\mathbb{B}_n(\rho)).
\end{eqnarray}
In particular, $(X,d,\mu)$ has the $n$-dimensional volume growth
\begin{eqnarray*}
C_{0}^{-1}\(C^{-1}\Phi\)^{\(\frac{\theta}{p}+\frac{1-\theta}{q}-\frac{1}{r}\)^{-1}}w_{n}\rho^{n}\leq
\mu(B(x_0,\rho))\leq C_{0}w_{n}\rho^{n}
\end{eqnarray*}
for all $\rho>0$, where $w_{n}$ is the volume of the unit ball in
$\mathbb{R}^{n}$.
\end{theorem}

\begin{remark}  \label{remark1.3}
\rm{ (1). When $q=\frac{(n-1)p}{n-p}$, then $\theta=1$,
$r=\frac{np}{n-p}$, and correspondingly, the Gagliardo-Nirenberg
inequality (\ref{a1}) degenerates into the following Sobolev
inequality
\begin{eqnarray*}
\left(\inx |u(x)|^{\frac{np}{n-p}}\dm\right)^{\frac{n-p}{np}}\leq
C\left(\inx |\n u|^p(x)\dm\right)^{\frac{1}{p}}
\end{eqnarray*}
for $u\in \mathrm{Lip_0}(X)$.
\\(2). The non-compactness of $(X,d)$ can be
assured by the validity of (\ref{a1}). In fact, if  $(X,d)$ is
bounded, then one can choose $q=\frac{(n-1)p}{n-p}$, then $\theta=1$
and $r=\frac{np}{n-p}$, which lets (\ref{a1}) become the Sobolev
inequality mentioned above, and in this setting, $u+\ell$ with
$\ell\rightarrow\infty$ clearly violates the validity of (\ref{a1}).
 $\\$ (3). By (\ref{a2}), we have
 \begin{eqnarray*}
\frac{\mu (B(x_0,R))}{w_{n}R^{n}}\leq C_0\frac{\mu
(B(x_0,r))}{w_{n}r^{n}}=C_0\frac{\mu (B(x_0,r))}{\mu_{E}
(\mathbb{B}_{n}(r))}
 \end{eqnarray*}
for $x_0\in{X}$ and $0<r<R$. Fixing $R$ and letting $r$ tends to
zero, by the volume condition (\ref{a3}) which describes the volume
behavior near $x_0$, we can obtain
\begin{eqnarray*}
\frac{\mu (B(x_0,R))}{w_{n}R^{n}}\leq C_0\cdot\liminf_{r\rightarrow
0}\frac{\mu (B(x_0,r))}{\mu_{E} (\mathbb{B}_{n}(r))}=1,
\end{eqnarray*}
which implies that $\mu (B(x_0,R))\leq C_0 w_{n}R^{n}$ for any
$R>0$. So, one can get the $n$-dimensional volume growth, i.e., the
last assertion of Theorem \ref{theorem1.1}, directly provided
(\ref{a4}) is proven.
 $\\$(4). As pointed out in \cite[Remark 1.2 (b)]{ko}, if
 $(X,d,\mu)$ satisfies the \emph{volume doubling condition}
 \begin{eqnarray*}
 \mu(B(x,2r))\leq\Lambda\mu(B(x,r)), \qquad
 \mathrm{for~some}~\Lambda\geq1~\mathrm{and~all}~x\in{X},~r>0,
 \end{eqnarray*}
then it is easy to get that the volume condition (\ref{a2}) is
satisfied with, e.g., $n\geq\log_{2}\Lambda$ and $C_{0}=1$.
Therefore, (\ref{a2}) can be comprehended as the volume doubling
condition with the explicit exponent $n$. Besides, one can regard
the volume condition (\ref{a3}) as a generalization of the classical
Bishop-Gromov volume comparison for complete manifolds with
non-negative Ricci curvature.
 $\\$(5). The assertion of having $n$-dimensional volume growth
 implies that, for instance, the cylinder
 $\mathbb{S}^{n-1}\times\mathbb{R}$ does not satisfy (\ref{a1}) for
 any $x\in X$ and $C$. The volume doubling condition (\ref{a2}) implies
 that the Hausdorff dimension $\mathrm{dim}_{H}X$ of $(X,d)$ is at
 most $n$. Besides, as in (3), by the volume conditions (\ref{a2}) and (\ref{a3}), we have
\begin{eqnarray*}
\frac{\mu (B(x_0,R))}{\mu_{E} (\mathbb{B}_{n}(R))} \leq C_0\frac{\mu
(B(x_0,r))}{\mu_{E} (\mathbb{B}_{n}(r))}
\end{eqnarray*}
for $x_0\in{X}$ and $0<r<R$, which implies that
\begin{eqnarray*}
\limsup_{R\rightarrow 0}\frac{\mu (B(x_0,R))}{\mu_{E}
(\mathbb{B}_{n}(R))}\leq\liminf_{r\rightarrow 0}C_0\frac{\mu
(B(x_0,r))}{\mu_{E} (\mathbb{B}_{n}(r))}=C_0.
\end{eqnarray*}
Therefore, we know that the \emph{Ahlfors $n$-regularity} at $x_0$
in the sense that $\eta^{-1}r^{n}\leq\mu (B(x_0,r))\leq\eta r^{n}$
for some $\eta\geq1$ and small $r>0$, which means that
$\mathrm{dim}_{H}X=n$. The volume doubling condition and the Ahlfors
regularity are important in analysis on metric measure spaces. For
this fact, see, e.g., \cite{hj} for the details. Note that the
choice of the constant $1$ chosen at the right hand side of
(\ref{a3}) is only for simplicity. In fact, by (\ref{a2}), we know
that $\eta_{x_0}:=\liminf_{r\rightarrow 0}\frac{\mu
(B(x_0,r))}{\mu_{E} (\mathbb{B}_{n}(r))}$ is positive. So, one can
normalize $\mu$ so as to satisfy (\ref{a3}) once $\eta_{x_0}$ is
bounded.
   }
\end{remark}

By Theorem \ref{BGF} (equivalently, see also Shen \cite{shz} or Ohta
\cite{os}), we know that for Finsler manifolds with non-negative
$n$-Ricci curvature (for this notion, see Definition \ref{defnr} for
the precise statement), the volume doubling condition (\ref{a2})
holds with $C_0=1$. For complete Finsler manifolds with non-negative
$n$-Ricci curvature, when the Gagliardo-Nirenberg inequality
(\ref{a1}) is satisfied with the best constant (i.e., $C=\Phi$), by
applying Theorems \ref{theorem1.2} and \ref{BGF}, we can prove the
following rigidity theorem.

\begin{corollary}\label{corollary1.4}
Let $(X,F)$ be a complete $n$-dimensional Finsler manifold. Let $p$,
$q$, $r$, $\theta$, $n$ be as in Theorem \ref{theorem1.1}, $x_0\in
X$, and $C_0\geq 1$.  Fix a positive smooth measure $\mu$ on $X$ and
assume that the $n$-Ricci curvature $\mathrm{Ric}_{n}$ of
$(X,F,\mu)$ is nonnegative. If the Gagliardo-Nirenberg inequality
(\ref{a1}) is satisfied with the best constant (i.e., $C=\Phi$) and
$\lim_{r\rightarrow 0}\frac{\mu (B(x_0,r))}{w_{n}r^{n}}=1$, then
under the volume doubling condition (\ref{a2}), we have the flag
curvature of $(X,F)$ is identically zero.
\end{corollary}

\begin{remark}\rm{
Finsler manifolds are special metric measure spaces with prescribed
Finsler structures. See Subsection 2.2 for a brief introduction to
Finsler manifolds.}
\end{remark}

A smooth metric measure space, which is also known as the weighted
measure space, is actually a Riemannian manifold equipped with some
measure (which is conformal to the usual Riemannian measure). More
precisely, for a given complete $n$-dimensional Riemannian manifold
$(M,g)$ with the metric $g$, the triple $(M,g,e^{-f}dv_{g})$ is
called a smooth metric measure space, with $f$ a \emph{smooth
real-valued} function on $M$ and $dv_{g}$ the Riemannian volume
element related to $g$ (sometimes, we also call $dv_{g}$ the volume
density). For a geodesic ball $B(x,r)$, we can define its
\emph{weighted (or $f$-)volume} $\mathrm{Vol}_f[B(x,r)]$ as follows
\begin{eqnarray}  \label{wv}
\mathrm{Vol}_f[B(x,r)]=\int_{B(x,r)}e^{-f}dv_{g}.
\end{eqnarray}
On a smooth metric measure space $(M,g,e^{-f}dv_{g})$, the so called
\emph{$\infty$-Bakry-\'{E}mery Ricci tensor} $\mathrm{Ric}_{f}$ is
defined by
\begin{eqnarray*}
\mathrm{Ric}_{f}=\mathrm{Ric}+\mathrm{Hess}f,
\end{eqnarray*}
which is also called the \emph{weighted Ricci curvature}. Bakry and
\'{E}mery \cite{be1,be2} introduced firstly and investigated
extensively the generalized Ricci tensor above and its relationship
with diffusion processes. In 2014, Mao \cite{m3} studied complete
non-compact smooth measure metric spaces with nonnegative weighted
Ricci curvature on which some Gagliardo-Nirenberg type inequality is
satisfied, and proved the following result.

 \begin{theorem} (\cite{m3}) \label{theorem1.5}
 Let $p$, $q$, $r$, $\theta$, $n$ be as in Theorem \ref{theorem1.1},
 and let $\(M,g,e^{-f}dv_{g}\)$ be an $n$-dimensional ($n\geq 2$) complete noncompact smooth metric measure space with non-negative
 weighted Ricci curvature. For a point $x_0\in M$ at which $f(x_0)$ is away from $-\infty$, assume that the radial
 derivative $\p_t f$ satisfies
$\p_t f\geq 0$ along all minimal geodesic segments from $x_0$, with
$t:=d(x_o,\cdot)$ the distance to $x_0$. Furthermore, for any $u\in
C_0^\infty(M)$ and some constant $C>0$, if the following
Gagliardo-Nirenberg type inequality
\begin{eqnarray}\label{a5}
\left(\inm |u(x)|^r e^{-f}\dn\right)^{\frac{1}{r}}\leq C\left(\inm
|\n u|^p(x)e^{-f}\dn\right)^{\frac{\theta}{p}}\left(\inm
|u(x)|^qe^{-f}\dn\right)^{\frac{1-\theta}{q}}
 \end{eqnarray}
 is satisfied, then we have
\begin{eqnarray}\label{a6}
\mathrm{Vol}_f[B(x_0,R)]\geq
e^{-f(x_0)}\(C^{-1}\Phi\)^{\(\frac{\theta}{p}+\frac{1-\theta}{q}-\frac{1}{r}\)^{-1}}V_0(R),\qquad
\forall R>0,
\end{eqnarray}
where $V_0(R)$ denotes the volume of an R-ball in $\R^n$.
\end{theorem}

\begin{remark}
\rm{ (1) By applying Theorem \ref{theorem1.2}, we can give an
\emph{alternative} proof to Theorem \ref{theorem1.5} for smooth
metric measure spaces of dimension $n\geq2$ -- see Subsection 3.3
for the details.
 $\\$(2). If the Gagliardo-Nirenberg type inequality (\ref{a5}) is
 satisfied with the best constant (i.e., $C=\Phi$), then by Theorems
 \ref{theorem1.5} and \ref{theorem3.1}, and together with
 generalized Bishop-type volume comparisons (cf. \cite[Theorem 3.3, Corollary 3.5 and Theorem
 4.2]{fmi}) for complete manifolds with \emph{radial} curvature
 bounded, a rigidity conclusion, $(M,g)$ is isometric
to $\left(\mathbb{R}^{n},g_{\mathbb{R}^{n}}\right)$ with
$g_{\mathbb{R}^{n}}$ being the usual Euclidean metric, can be
obtained (cf. \cite[Corollary 1.5]{m3} for the precise statement).

 }
\end{remark}

It is interesting to know \emph{under what kind of conditions} a
complete open $n$-manifold ($n\geqslant2$) is isometric to
$\mathbb{R}^{n}$ or has finite topological type, which in essence
has relation with the splittingness of the prescribed manifold. This
is a classical topic in the global differential geometry, which has
been investigated intensively (see, e.g., \cite{dcx2,m,pp}).

\section{Proofs of Theorem \ref{theorem1.2} and Corollary \ref{corollary1.4}}
\renewcommand{\thesection}{\arabic{section}}
\renewcommand{\theequation}{\thesection.\arabic{equation}}
\setcounter{equation}{0}

\subsection{Proof of Theorem \ref{theorem1.2}}

{\emph{Proof}}. As pointed out in Remark \ref{remark1.3} (3), if we
want to get the $n$-dimensional volume growth assertion in Theorem
\ref{theorem1.2}, we only need to show (\ref{a4}). Now, in the rest
of this subsection, we would like to give the details of the proof
of (\ref{a4}) as follows.

 First, we introduce two auxiliary
functions $F,G: (0,+\infty)\rightarrow \mathbb{R}$ defined by
\begin{eqnarray*}
F(\lambda):=\inx
\frac{1}{\(\la+d(x_0,x)^{\frac{p}{p-1}}\)^\frac{(p-1)q}{q-p}}\dm
\end{eqnarray*}
and
\begin{eqnarray*}
 G(\lambda):=\inr
\frac{1}{\(\la+|x|^{\frac{p}{p-1}}\)^\frac{(p-1)q}{q-p}}\de
\end{eqnarray*}
respectively, which are well defined and of class $C^1$.

By the layer cake representation of functions, one has
\begin{eqnarray*}
F(\lambda)=\ino \mu\left\{x\in
X:\frac{1}{\(\la+d(x_0,x)^{\frac{p}{p-1}}\)^\frac{(p-1)q}{q-p}}>s\right\}ds.
\end{eqnarray*}
By taking into account that $\mathrm{diam} (X)=\infty$ and making
the variable change
$$s=\frac{1}{\(\la+\rho^{\frac{p}{p-1}}\)^\frac{(p-1)q}{q-p}},$$
then
\begin{eqnarray}\label{b1}
F(\lambda)=\frac{pq}{q-p}\ino\mu (B(x_0,\rho))f(\la,\rho)d\rho,
\end{eqnarray}
where
$$f(\la,\rho)=\frac{\rho^{\frac{1}{p-1}}}{\(\la+\rho^{\frac{p}{p-1}}\)^\frac{(q-1)p}{q-p}}.$$
Similar to the above process, we can also get
\begin{eqnarray}\label{b2}
G(\lambda)=\frac{pq}{q-p}\ino\mu_E
(\mathbb{B}_n(\rho))f(\la,\rho)d\rho.
\end{eqnarray}

On the other hand, since the inequalities (\ref{a2}) and (\ref{a3})
hold on $(X,d,\mu)$, we have
\begin{eqnarray}\label{b3}
\mu(B(x_0,\rho))\leq C_0 \mu_E(\mathbb{B}_n(\rho)).
\end{eqnarray}
Then, it follows from (\ref{b1})--(\ref{b3}) that
\begin{eqnarray}\label{b4}
0\leq F(\lambda)\leq C_0 G(\la).
\end{eqnarray}
Since $q<\frac{np}{n-p}$,  we have
$$n+\frac{1}{p-1}-\frac{p^2(q-1)}{(p-1)(q-p)}=n+\frac{1}{p-1}-\frac{p^2}{p-1}-\frac{p^2}{q-p}<-1,$$
and from which we know that $0\leq F(\la)< \infty$, for any
$\lambda> 0,$ and $F(\la)$ is differentiable. Also, we have
\begin{eqnarray}\label{b41}
F'(\la)=-\frac{(p-1)q}{q-p}\inx
\frac{\dm}{\(\la+d(x_0,x)^{\frac{p}{p-1}}\)^{\frac{(q-1)p}{q-p}}}.
\end{eqnarray}

For each $\la>0$, consider the sequence of functions
\begin{eqnarray*}
u_{\la,k}(x):=\mathrm{max}\{0,\mathrm{min}\{0,k-d(x_0,x)\}+1\}\(\la+\mathrm{max}\{d(x_0,x),k^{-1}\}^{\frac{p}{p-1}}\)^{-\frac{p-1}{q-p}}.
\end{eqnarray*}
Since $(X,d)$ is proper, $\mathrm{supp}(u_{\la,k}):=\{x\in
X:d(x_0,x)\leq k+1\}$ is compact. Therefore, we have $u_{\la,k}\in
\mathrm{Lip}_0(X)$ for every $\la>0$ and $k\in \mathbb{N}$. Set
\begin{eqnarray*}
u_\la(x):=\lim_{k\rightarrow\infty}u_{\la,k}(x)=\(\la+d(x_0,x)^{\frac{p}{p-1}}\)^{-\frac{p-1}{q-p}}.
\end{eqnarray*}
Since the functions $u_{\la,k}$ verify the Gagliardo-Nirenberg
inequality (\ref{a1}), by an approximation based on (\ref{b4}), we
know that $u_\la$ verifies the Gagliardo-Nirenberg inequality
(\ref{a1}) also. Together with the fact that $x\rightarrow d(x_0,x)$
is 1-Lipschitz (i.e., $|Dd(x_0,\cdot)|(x)\leq1$ for all $x$), we can
obtain
\begin{eqnarray}\label{b5}
\nonumber&&\left(\inx
\frac{\dm}{\(\la+d(x_0,x)^{\frac{p}{p-1}}\)^{\frac{(q-1)p}{q-p}}}\right)^{\frac{1}{r}}\\&\leq&
 C\(\frac{p}{q-p}\)^\theta\left(\inx \frac{d(x_0,x)^{\frac{p}{p-1}}\dm}{\(\la+d(x_0,x)^{\frac{p}{p-1}}\)^{\frac{(q-1)p}{q-p}}}\right)^{\frac{\theta}{p}}\left(\inx
 \frac{\dm}{\(\la+d(x_0,x)^{\frac{p}{p-1}}\)^{\frac{(q-1)p}{q-p}}}\right)^{\frac{1-\theta}{q}}
 \qquad
 \end{eqnarray}
by using a chain rule for the local Lipschitz constant. By the
definition of $F(\la)$ and (\ref{b41}), the above equality can be
rewritten as follows
\begin{eqnarray}\label{b6}
\(-F'(\la)\)^{\frac{p}{\theta r}}\leq
l\(F(\la)+\frac{q-p}{(p-1)q}\la
F'(\la)\)F(\la)^{\frac{(1-\theta)p}{\theta q}},
 \end{eqnarray}
where\begin{eqnarray*}
l=C^{\frac{p}{\theta}}\(\frac{p}{q-p}\)^p\(\frac{(p-1)q}{q-p}\)^{\frac{p}{\theta
r}}.
\end{eqnarray*}

Since $v_\la(x)=\(\la+|x|^{\frac{p}{q-p}}\)^{-\frac{p-1}{q-p}}$ is a
minimizer of the Gagliardo-Nirenberg inequality in $\mathbb{R}^n$,
then for every $\la>0$, the following equality
\begin{eqnarray*}
\left(\inr |v_\la(x)|^r\de\right)^{\frac{1}{r}}= \Phi\left(\inr |\n
v_\la|^p(x)\de\right)^{\frac{\theta}{p}}\left(\inr
|v_\la(x)|^q\de\right)^{\frac{1-\theta}{q}}
 \end{eqnarray*}
holds. By the definition of $G(\la)$ and a similar argument as
(\ref{b6}), the above equality can be rewritten as follows
\begin{eqnarray}\label{b8}
\(-G'(\la)\)^{\frac{p}{\theta
r}}=\widetilde{l}\(G(\la)+\frac{q-p}{(p-1)q}\la
G'(\la)\)G(\la)^{\frac{(1-\theta)p}{\theta q}},
 \end{eqnarray}
where\begin{eqnarray*}
\widetilde{l}=\Phi^{\frac{p}{\theta}}\(\frac{p}{q-p}\)^p\(\frac{(p-1)q}{q-p}\)^{\frac{p}{\theta
r}}.
\end{eqnarray*}
Substituting $G(\la)=G(1)\la^{(p-1)\(\frac{n}{p}-\frac{q}{q-p}\)}$
into (\ref{b8}), we have
\begin{eqnarray}\label{b9}
\(1-\frac{n(q-p)}{pq}\)^{\frac{p}{\theta
r}}=\Phi^{\frac{p}{\theta}}\(\frac{p}{q-p}\)^p\(\frac{(q-p)^n}{pq}\)G(1)^{\frac{p}{\theta}\(\frac{\theta}{p}+\frac{1-\theta}{q}-\frac{1}{r}\)}.
\end{eqnarray}
Consider the constant $A$ given by
\begin{eqnarray}\label{b10}
\(1-\frac{n(q-p)}{pq}\)^{\frac{p}{\theta
r}}=C^{\frac{p}{\theta}}\(\frac{p}{q-p}\)^p\(\frac{(q-p)^n}{pq}\)A^{\frac{p}{\theta}\(\frac{\theta}{p}+\frac{1-\theta}{q}-\frac{1}{r}\)}.
\end{eqnarray}
It is easy to check that the function
\begin{eqnarray*}
H_0(\la)=A\la^{(p-1)\(\frac{n}{p}-\frac{q}{q-p}\)},\qquad \la\in
(0,+\infty),
\end{eqnarray*}
satisfies the differential equation
\begin{eqnarray}\label{b11}
\(-q'(\la)\)^{\frac{p}{\theta r}}=l\(q(\la)+\frac{q-p}{(p-1)q}\la
q'(\la)\)q(\la)^{\frac{(1-\theta)p}{\theta q}}.
 \end{eqnarray}
By (\ref{b9}) and (\ref{b10}), we get
\begin{eqnarray*} A=\(\frac{\Phi}{C}\)^{\(\frac{\theta}{p}+\frac{1-\theta}{q}-\frac{1}{r}\)^{-1}}G(1),
\end{eqnarray*}
which implies
\begin{eqnarray}\label{b12}
H_0(\la)=\(\frac{\Phi}{C}\)^{\(\frac{\theta}{p}+\frac{1-\theta}{q}-\frac{1}{r}\)^{-1}}\la^{(p-1)\(\frac{n}{p}-\frac{q}{q-p}\)}G(1)
=\(\frac{\Phi}{C}\)^{\(\frac{\theta}{p}+\frac{1-\theta}{q}-\frac{1}{r}\)^{-1}}G(\la).
\end{eqnarray}

In the following, we will show that when $C>\Phi$, for every
$\la>0$,
\begin{eqnarray}\label{b13}
F(\la)\geq H_0(\la).
\end{eqnarray}
First, we \emph{claim}  that if $F(\la_0)< H_0(\la_0)$, for some
$\la_0>0$, then $F (\la)<H_0(\la), \forall \la \in (0, \la_0].$ We
prove this by contradiction. Suppose that the claim is not true.
Then there exists some $\widetilde{\la}\in (0,\la_0)$ such that $F
(\widetilde{\la})\geq H_0(\widetilde{\la})$. Set
$\la_1:=\mathrm{sup}\{\la<\la_0:F(\la)=H_0(\la)\}$. Then for any
$\la\in [\la_1,\la_0], 0<F(\la)\leq H_0(\la)$, and so from
(\ref{b6}), we have
\begin{eqnarray}\label{b14}
\(-F'(\la)\)^{\frac{p}{\theta r}}\leq
l\(H_0(\la)+\frac{q-p}{(p-1)q}\la
F'(\la)\)H_0(\la)^{\frac{(1-\theta)p}{\theta q}}.
 \end{eqnarray}
For every $\la>0$, define a function $z_\la:(0,\infty)\rightarrow
\mathbb{R}$ by $z_\la(\rho)=\rho^{\frac{p}{\theta
r}}+\frac{l\la(q-p)\rho}{(p-1)q}H_0(\la)^{\frac{(1-\theta)p}{\theta
q}}$. Clearly, $z_\lambda$ is increasing. Hence, when $\la\in
[\la_1,\la_0],$ we infer from (\ref{b14}) and (\ref{b11}) that
\begin{eqnarray*}
z_\la(-F'(\la))&=&(-F'(\la))^{\frac{p}{\theta r}}+\frac{l\la(q-p)}{(p-1)q}(-F'(\la))H_0(\la)^{\frac{(1-\theta)p}{\theta q}}\\
&\leq&lH_0(\la)^{1+\frac{(1-\theta)p}{\theta q}}=z_\la(-H'_0(\la)),
 \end{eqnarray*}
which means $F'(\la)\geq H_{0}'(\la), \forall \la\in [\la_1,\la_0].$
Thus, we know that the function $F-H_0$ is increasing on
$[\la_1,\la_0]$, which implies that
\begin{eqnarray*}
0=(F-H_0)(\la_1)\leq(F-H_0)(\la_0)<0.
 \end{eqnarray*}
This is a contradiction. Hence, the above \emph{claim} is true.

By (\ref{a3}), we know that for every $\varepsilon>0$, there exists
$\delta>0$ such that $\mu(B(x_0,\rho))\geq
(1-\varepsilon)\mu_E(\mathbb{B}(\rho))$ for all $0\leq\rho\leq
\delta.$ Therefore, by (\ref{b6}) and making a variable change
$\rho=\la^{\frac{1}{2-ap}}t$, we can get
\begin{eqnarray*}
F(\la)&\geq&\frac{pq}{p-q}(1-\varepsilon)\int_0^\delta
\mu_E(\mathbb{B}_n(\rho)f(\la,\rho))d\rho\\&=&
\frac{1-a}{n-1+a}(1-\varepsilon)\la^{\frac{(p-1)n}{p}+1-\frac{p(q-1)}{q-p}}\int_0^{\delta\la^{\frac{1-p}{p}}}\mu_E(\mathbb{B}_n(t)f(1,t))dt.
 \end{eqnarray*}
On the other hand, we have
\begin{eqnarray*}
G(\la)=\frac{1-a}{n-1+a}\la^{\frac{(p-1)n}{p}+1-\frac{p(q-1)}{q-p}}\int_0^\infty\mu_E(\mathbb{B}_n(t)f(1,t))dt.
\end{eqnarray*}
Therefore, we have
\begin{eqnarray*}
\liminf_{\la\rightarrow 0}\frac{F(\la)}{G(\la)}\geq 1-\varepsilon.
\end{eqnarray*}
Letting $\varepsilon\rightarrow0$ yields
\begin{eqnarray*}
\liminf_{\la\rightarrow 0}\frac{F(\la)}{G(\la)}\geq 1.
\end{eqnarray*}
When $C>\Phi$, we infer from the above inequality and (\ref{b12})
that
\begin{eqnarray*}
\liminf_{\la\rightarrow
0}\frac{F(\la)}{H_0(\la)}=\(\frac{C}{\Phi}\)^{\(\frac{\theta}{p}+\frac{1-\theta}{q}-\frac{1}{r}\)^{-1}}\liminf_{\la\rightarrow
0}\frac{F(\la)}{G(\la)}\geq
\(\frac{C}{\Phi}\)^{\(\frac{\theta}{p}+\frac{1-\theta}{q}-\frac{1}{r}\)^{-1}}>
1.
\end{eqnarray*}
Then, together with the previous \emph{claim}, we can get
$F(\la)\geq H_0(\la), \forall \la>0$.  Thus, for any $\la>0$, we can
get from (\ref{b1}), (\ref{b2}), (\ref{b12}) that
\begin{eqnarray}\label{b15}
\ino\left
\{\mu(B(x_0,\rho))-b\mu_E(\mathbb{B}_n(\rho))\right\}f(\la,\rho)d\rho\geq
0,
\end{eqnarray}
where
$b=\(C^{-1}\Phi\)^{\(\frac{\theta}{p}+\frac{1-\theta}{q}-\frac{1}{r}\)^{-1}}$.
By (\ref{a2}), for a fixed $\rho>0$, we have
\begin{eqnarray*}
C_0\frac{\mu(B(x_0,\rho))}{\mu_E(\mathbb{B}_n(\rho))}\geq
\frac{\mu(B(x_0,r))}{\mu_E(\mathbb{B}_n(r))}
\end{eqnarray*}
for any $r>\rho\geq0$. We can assume
\begin{eqnarray*}
b_0:=\limsup_{r\rightarrow\infty}\frac{\mu(B(x_0,r))}{\mu_E(\mathbb{B}_n(r))}.
\end{eqnarray*}
In order to prove (\ref{a4}) in the case that $C>\Phi$, it suffices
to show that $b_0\geq b$. We will prove this by contradiction. By
the definition of $b_0$, we know that for some $\rho_0>0$, there
exists $\varepsilon_0>0$ such that
\begin{eqnarray}\label{b16}
\frac{\mu(B(x_0,\rho))}{\mu_E(\mathbb{B}_n(\rho))}\leq
b-\varepsilon_0, \qquad \forall \rho\geq \rho_0.
\end{eqnarray}
Substituting (\ref{b16}) into (\ref{b15}), and together with
(\ref{b3}), for every $\la>0$, we have
\begin{eqnarray*}
0&\leq&\ino \left\{\mu(B(x_0,\rho))-b\mu_E(\mathbb{B}_n(\rho))\right\}f(\la,\rho)d\rho\\
&\leq&\int_0^{\rho_0}\mu(B(x_0,\rho))f(\la,\rho)d\rho
+( b-\varepsilon_0)\int_{\rho_0}^{+\infty} \mu_E(\mathbb{B}_n(\rho))f(\la,\rho)d\rho\\&&-b\ino\mu_E(\mathbb{B}_n(\rho))f(\la,\rho)d\rho\\
&\leq& C_0\int_0^{\rho_0} \mu_E(\mathbb{B}_n(\rho))f(\la,\rho)d\rho-b\int_0^{\rho_0}\mu_E(\mathbb{B}_n(\rho))f(\la,\rho)d\rho\\
&&-\varepsilon_0\int_{\rho_0}^{+\infty}
\mu_E(\mathbb{B}_n(\rho))f(\la,\rho)d\rho\\&=&
\(C_0-b+\varepsilon_0\)\int_0^{\rho_0}
\mu_E(\mathbb{B}_n(\rho))f(\la,\rho)d\rho-\varepsilon_0\int_{0}^{+\infty}
\mu_E(\mathbb{B}_n(\rho))f(\la,\rho)d\rho
\\&=&
\(C_0-b+\varepsilon_0\)\int_0^{\rho_0}
\mu_E(\mathbb{B}_n(\rho))f(\la,\rho)d\rho-\varepsilon_0\frac{q-p}{pq}\la^{(p-1)\(\frac{n}{p}-\frac{q}{q-p}\)}G(1).
\end{eqnarray*}
Since
$f(\la,\rho)=\frac{\rho^{\frac{1}{p-1}}}{\(\la+\rho^{\frac{p}{p-1}}\)^\frac{(q-1)p}{q-p}},$
one has
\begin{eqnarray*}
\int_0^{\rho_0} \rho^n
f(\la,\rho)d\rho=\int_0^{\rho_0}\frac{\rho^{n+\frac{1}{p-1}}}{\(\la+\rho^{\frac{p}{p-1}}\)^\frac{(q-1)p}{q-p}}d\rho\leq
\la^{\frac{-p(q-1)}{q-p}}\frac{\rho_0^{n+1+\frac{1}{p-1}}}{\(n+1+\frac{1}{p-1}\)}.
\end{eqnarray*}
By the above two inequalities, we can get the inequality of type
\begin{eqnarray*}
M_1\la^{(p-1)\(\frac{n}{p}-\frac{q}{q-p}\)}\leq M_2
\la^{\frac{-p(q-1)}{q-p}}, \qquad \forall \la >0,
\end{eqnarray*}
where $M_1, M_2>0$ are constants independent of $\lambda$. Observing
$\frac{-p(q-1)}{q-p}-(p-1)\(\frac{n}{p}-\frac{q}{q-p}\)=\frac{(1-p)n}{p}-1<0$,
letting $\la\rightarrow +\infty$ in the above inequality, one can
obtain a contradiction. This means that (\ref{a4}) holds in the case
that $C>\Phi$.

When $C=\Phi$, we can also get (\ref{a4}). In fact, in this case,
for any fixed $\delta>0$, we have
\begin{eqnarray*}
\left(\inx |u(x)|^r\dm\right)^{\frac{1}{r}}\leq
(\Phi+\delta)\left(\inx |\n
u|^p(x)\dm\right)^{\frac{\theta}{p}}\left(\inx
|u(x)|^q\dm\right)^{\frac{1-\theta}{q}}.
 \end{eqnarray*}
 Therefore, for any $x\in X$, by the previous argument, we have
\begin{eqnarray*}
 \mu (B(x,\rho))\geq C_{0}^{-1}
 \(\frac{\Phi}{\Phi+\delta}\)^{\(\frac{\theta}{p}+\frac{1-\theta}{q}-\frac{1}{r}\)^{-1}}\mu_E(\mathbb{B}_n(\rho)),
 \qquad\forall\rho>0.
\end{eqnarray*}
Letting $\delta\rightarrow0$, we can obtain
\begin{eqnarray*}
\mu (B(x,\rho))\geq C_{0}^{-1}\mu_E(\mathbb{B}_n(\rho)),
 \qquad\forall\rho>0,
\end{eqnarray*}
which implies (\ref{a4}) holds in the case that $C=\Phi$.

This completes the proof of Theorem \ref{theorem1.2}.  \hfill $\Box$

\subsection{Preliminary notions and a Bishop-Gromov type volume comparison theorem in Finsler geometry}

Before applying Theorem \ref{theorem1.2} to prove Corollary
\ref{corollary1.4}, we briefly recall some concepts in Finsler
geometry. We refer to \cite{BCS} for a fundamental but overall
introduction about Finsler geometry.

\begin{defn} \label{defnr} \rm{
Let $X$ be a connected $n$-dimensional smooth manifold and
$TX=\bigcup_{x\in X}T_{x}X$ be its tangent bundle. The pair $(X,F)$
is called a \emph{Finsler manifold} if a continuous function
$F:TX\rightarrow[0,\infty)$ satisfies the following conditions
 $\\$(1) $F\in C^{\infty}\(TX\backslash\{0\}\)$;
 $\\$(2) $F(x,tv)=|t|F(x,v)$ for all $t\in\mathbb{R}$ and $(x,v)\in
 TX$;
 $\\$(3) The $n\times n$ matrix
 \begin{eqnarray} \label{2-2-1}
 g_{ij}(x,v):=\frac{1}{2}\frac{\partial^{2}\(F^2\)}{\partial v^{i}\partial
 v^{j}}, \qquad \mathrm{where}~v=\sum\limits_{i=1}^{n}v^{i}\frac{\partial}{\partial
 x^{i}},
 \end{eqnarray}
 is positive definite for all $(x,v)\in  TX\backslash\{0\}$. $F$ is
 called the Finsler structure of $(X,F)$.
}
\end{defn}

We will denote by $\langle,\rangle_v$ the inner product on $T_{x}X$
induced by (\ref{2-2-1}). We know that $(X,F)$ becomes a Riemannian
manifold if and only if $g_{ij}(x,v)$ is independent of $v$ in each
$T_{x}X\backslash\{0\}$. For a smooth curve $\sigma:[0,l)\rightarrow
X$, one can define its \emph{integral length} $L_{F}\sigma$ by
$L_{F}\sigma=\int_{0}^{l}F(\sigma,\dot{\sigma})dt$. Based on this,
the \emph{distance function} $d_{F}:X\times X\rightarrow[0,\infty)$
can be defined by $d_{F}(x_1,x_2):=\inf_{\sigma}L_{F}\sigma$, where
$\sigma$ runs over all smooth curves from $x_1$ to $x_2$. A smooth
curve $\sigma:[0,l)\rightarrow X$ is called a \emph{geodesic} if it
locally minimizes $d_F$ and has a constant speed (i.e.,
$F(\sigma,\dot{\sigma})$ is constant). The geodesic (Euler-Lagrange)
equation can be written down in terms of covariant derivative along
$\sigma$ (see \cite{BCS} for the details). The Finsler manifold
$(X,F)$ is complete if any geodesic $\sigma:[0,l)\rightarrow X$ can
be extended to a geodesic $\sigma:\mathbb{R}\rightarrow X$.

Like the Riemannian case, we can also do the geodesic variation in
the Finsler case. In fact, let
$\sigma:(-\epsilon,\epsilon)\times[0,l]\rightarrow X$ be a smooth
geodesic variation (i.e., $t\rightarrow\sigma(s,t)$ is geodesic for
each $s$), and set $\eta(t)=\sigma(0,s)$. Then the variational
vector field $J(t):=\frac{\partial\sigma}{\partial s}(0,t)$
satisfies the following Jacobi equation
\begin{eqnarray} \label{2-2-2}
D_{\dot{\eta}}^{\dot{\eta}}D_{\dot{\eta}}^{\dot{\eta}}J+R^{\dot{\eta}}(J,\dot{\eta})\dot{\eta}=0,
\end{eqnarray}
where $D^{\dot{\eta}}$ is the covariant derivative w.r.t. the vector
$\dot{\eta}$, and $R^{\dot{\eta}}$ is the curvature tensor (see
\cite{BCS} for the details). For vectors $v,w\in T_{x}X$, which are
linearly independent, and $\mathcal {S}=\mathrm{span}\{v,w\}$, the
\emph{flag curvature} of the \emph{flag} $(\mathcal{S};v)$ can be
defined as follows
\begin{eqnarray*}
K(\mathcal{S},v):=\frac{\langle
R^{v}(w,v)v,w\rangle_{v}}{F(v)^{2}\langle w,w\rangle_{v}-\langle
v,w\rangle_{v}^{2}}.
\end{eqnarray*}
If $(X,F)$ is a Riemannian manifold, then the flag curvature
degenerates into the sectional curvature which only depends on
$\mathcal{S}$ (not on the choice of $v\in\mathcal{S}$). Choose $v\in
T_{x}X$ with $F(x,v)=1$, and let $\{e_{i}\}_{i=1}^{n}$ with
$e_{n}=v$ be an orthonormal basis of $(T_{x}X,\langle,\rangle_{v})$
with $\langle,\rangle_{v}$ induced from (\ref{2-2-1}). Set
$\mathcal{S}_{i}=\mathrm{span}\{e_{i},v\}$ for $i=1,2,\ldots,n-1$.
The Ricci curvature of $v$ is defined by
$\mathrm{Ric}(v):=\sum_{i=1}^{n-1}K(\mathcal{S}_{i};v)$. We also set
$\mathrm{Ric}(cv):=c^{2}\mathrm{Ric}(v)$ for $c\geq0$.

For those Finsler curvatures mentioned above, Shen has explained
them from the Riemannian viewpoint (see \cite[Section 6.2 of Chapter
6]{shz1}). Fixing $v\in T_{x}X\backslash\{0\}$ and extending it to a
smooth vector field $V$ around $x$ such that all integral curves of
$V$ are geodesics, then the flag curvature $K(\mathcal{S},v)$ is the
same as the sectional curvature of $\mathcal{S}$ w.r.t. the
Riemannian structure $\langle,\rangle_{V}$, and correspondingly,
$\mathrm{Ric}(v)$ is the same as the Ricci curvature of $v$ w.r.t.
$\langle,\rangle_{V}$. This fact leads to the following definition
of $N$-Ricci curvature associated with an arbitrary measure on $X$
(see also, e.g., \cite{ko,os} for this notion).

\begin{defn}  \rm{
Let $\mu$ be a positive smooth measure on $X$. Given $v\in
T_{x}X\backslash\{0\}$, let $\sigma:(-\epsilon,\epsilon)\rightarrow
X$ be the geodesic with $\dot{\sigma}=v$ and decompose $\mu$ along
$\sigma$ as $\mu=e^{-\psi}\mathrm{Vol}_{\dot{\sigma}}$, where
$\mathrm{Vol}_{\dot{\sigma}}$ is the volume element of the
Riemannian structure $\langle,\rangle_{\dot{\sigma}}$. Then, for
$N\in[n,\infty]$, the \emph{N-Ricci curvature} $\mathrm{Ric}_{N}$ is
defined by
\begin{eqnarray*}
\mathrm{Ric}_{N}(v)=\mathrm{Ric}(v)+(\psi\circ\sigma)''(0)-\frac{(\psi\circ\sigma)'(0)^{2}}{N-n},
\end{eqnarray*}
where the third term is understood as $0$ if $N=\infty$ or if $N=n$
with $(\psi\circ\sigma)'(0)=0$, and as $-\infty$ if $N=n$ with
$(\psi\circ\sigma)'(0)\neq0$.}
\end{defn}

By applying the concept of the $N$-Ricci curvature
$\mathrm{Ric}_{N}$, Ohta \cite{os} proved the following
Bishop-Gromov type volume comparison result in the Finsler case.

\begin{theorem} \label{BGF}
(\cite[Theorem 7.3]{os}) Let $(X,F,\mu)$ be a complete
$n$-dimensional Finsler manifold with nonnegative N-Ricci curvature.
Then we have
\begin{eqnarray*}
\frac{\mu (B(x,R))}{\mu (B(x,r))}\leq \left(\frac{R}{r}\right)^{N},
\quad \mathrm{for~every}~ x\in X, ~\mathrm{and}~ 0<r<R.
\end{eqnarray*}
Moreover, if equality holds with $N=n$ for all $x\in X$ and $0<r<R$,
then any Jacobi field $J$ along a geodesic $\sigma$ has the form
$J(t)=tP(t)$, where $P$ is a parallel vector field along $\sigma$
(i.e., $D_{\dot{\sigma}}^{\dot{\sigma}}P\equiv0$).
\end{theorem}

\subsection{Proof of Corollary \ref{corollary1.4}}

{\emph{Proof}}. Since $(X,F)$ is complete, by applying the
Hopf-Rinow theorem, we know that $(X,d_F,\mu)$ is a \emph{proper}
metric measure space. Since the $n$-Ricci curvature
$\mathrm{Ric}_{n}$ is nonnegative, by Theorem \ref{BGF}, we can
obtain (\ref{a2}) with $C_0=1$. As pointed out in Remark
\ref{remark1.3} (5), one can normalize the fixed positive measure
$\mu$ such that (\ref{a3}) is satisfied. Then by these two facts,
similar to Remark \ref{remark1.3} (3), we can easily get
\begin{eqnarray*}
\mu (B(x,\rho))\leq \mu_E(\mathbb{B}_n(\rho)), \qquad
\mathrm{for~all}~\rho>0,~x\in X.
\end{eqnarray*}
However, since the Gagliardo-Nirenberg inequality (\ref{a1}) is
satisfied with the best constant (i.e., $C=\Phi$), by Theorem
\ref{theorem1.2}, we have
\begin{eqnarray*}
\mu (B(x,\rho))\geq \mu_E(\mathbb{B}_n(\rho)), \qquad
\mathrm{for~all}~\rho>0,~x\in X.
\end{eqnarray*}
Therefore, $\mu (B(x,\rho))= \mu_E(\mathbb{B}_n(\rho))$ for all
$\rho>0$ and $x\in X$. By applying Theorem \ref{BGF} directly, we
know that every Jacobi field $J$ along a geodesic $\sigma$ has the
form $J(t)=tP(t)$, where $P$ is a parallel vector field along
$\sigma$. Together with the Jacobi equation (\ref{2-2-2}), it
follows that $R^{\dot{\sigma}}(J,\dot{\sigma})\dot{\sigma}\equiv0$.
Then $K(\mathcal{S};\dot{\sigma})\equiv0$ with
$\mathcal{S}=\mathrm{span}(\dot{\sigma},P)$. Since $\sigma$ and $J$
are arbitrary, we know $K\equiv0$, which equivalently says that the
flag curvature of $(X,F)$ is identically zero. \hfill $\Box$

\section{Proof of Theorem \ref{theorem1.5} }
\renewcommand{\thesection}{\arabic{section}}
\renewcommand{\theequation}{\thesection.\arabic{equation}}
\setcounter{equation}{0}  \label{intro}

In this section, as mentioned before, we would like to give an
alternative proof to Theorem \ref{theorem1.5}. However, before that,
we need to introduce some notions. For more details, we refer
readers to \cite{fmi,m1,m2,m3,m4}.

\subsection{Some basic notions}
Denote by $\mathbb{S}^{n-1}$ the unit sphere in $\mathbb{R}^{n}$.
Given an $n$-dimensional ($n\geqslant2$) complete Riemannian
manifold $(M,g)$ with the metric $g$, for a point $x\in{M}$, let
$S_{x}^{n-1}$ be the unit sphere with center $x$ in the tangent
space $T_{x}M$, and let $Cut(x)$ be the cut-locus of $x$, which is a
closed set of zero $n$-Hausdorff measure. Clearly,
\begin{eqnarray*}
\mathbb{D}_{x}=\left\{t\xi|0\leqslant{t}<d_{\xi},\xi\in{S_{x}^{n-1}}\right\}
\end{eqnarray*}
is a star-shaped open set of $T_{x}M$, and through which the
exponential map
$\exp_{x}:\mathbb{D}_{x}\rightarrow{M}\backslash{Cut(x)}$ gives a
diffeomorphism from $\mathbb{D}_{x}$ to the open set
$M\backslash{Cut(x)}$, where $d_{\xi}$ is defined by
\begin{eqnarray*}
d_{\xi}=d_{\xi}(x) : = \sup\{t>0|~\gamma_{\xi}(s):=
\exp_x(s\xi)~{\rm{is~ the~ unique}}~{\rm{minimal
~geodesic~joining}}~x ~ {\rm{and}}~\gamma_{\xi}(t)\}.
\end{eqnarray*}

We can introduce two
 important maps used to construct the geodesic spherical coordinate chart at a prescribed point
 on a Riemannian manifold. For a fixed vector $\xi\in{T_{x}M}$,
$|\xi|=1$, let $\xi^{\bot}$ be the orthogonal complement of
$\{\mathbb{R}\xi\}$ in $T_{x}M$, and let
$\tau_{t}:T_{x}M\rightarrow{T_{\exp_{x}(t\xi)}M}$ be the parallel
translation along $\gamma_{\xi}(t)$. The path of linear
transformations
$\mathbb{A}(t,\xi):\xi^{\bot}\rightarrow{\xi^{\bot}}$ is defined by
 \begin{eqnarray*}
\mathbb{A}(t,\xi)\eta=(\tau_{t})^{-1}Y_{\eta}(t),
 \end{eqnarray*}
where $Y_{\eta}(t)=d(\exp_x)_{(t\xi)}(t\eta)$ is the Jacobi field
along $\gamma_{\xi}(t)$ satisfying $Y_{\eta}(0)=0$, and
$(\nabla_{t}Y_{\eta})(0)=\eta$. Moreover, for $\eta\in{\xi^{\bot}}$,
set $
\mathcal{R}(t)\eta=(\tau_{t})^{-1}R(\gamma'_{\xi}(t),\tau_{t}\eta)
\gamma'_{\xi}(t), $ where the curvature tensor $R(X,Y)Z$ is defined
by $R(X,Y)Z=-[\nabla_{X},$ $ \nabla_{Y}]Z+ \nabla_{[X,Y]}Z$. Then
$\mathcal{R}(t)$ is a self-adjoint operator on $\xi^{\bot}$, whose
trace is the radial Ricci tensor
$\mathrm{Ric}_{\gamma_{\xi}(t)}\left(\gamma'_{\xi}(t),\gamma'_{\xi}(t)\right)$.
 Clearly, the map $\mathbb{A}(t,\xi)$ satisfies the Jacobi
equation
 $\mathbb{A}''+\mathcal{R}\mathbb{A}=0$ with initial conditions
 $\mathbb{A}(0,\xi)=0$, $\mathbb{A}'(0,\xi)=I$.
By Gauss's lemma, the Riemannian metric of $M \backslash Cut(x)$ in
the geodesic spherical coordinate chart can be expressed by
 \begin{eqnarray} \label{2.1}
 ds^{2}(\exp_{x}(t\xi))=dt^{2}+|\mathbb{A}(t,\xi)d\xi|^{2}, \qquad
\forall t\xi\in\mathbb{D}_{x}.
 \end{eqnarray}
We consider the metric components  $g_{ij}(t,\xi)$, $i,j\geq 1$, in
a coordinate system $\{t, \xi_a\}$ formed by fixing  an orthonormal
basis $\{\eta_a, a\geq 2\}$ of
 $\xi^{\bot}=T_{\xi}S^{n-1}_x$, and then extending it to a local frame $\{\xi_a, a\geq2\}$ of
$S_x^{n-1}$. Define
 a function $J>0$ on $\mathbb{D}_{x}\backslash\{x\}$ by
\begin{equation} \label{2.2}
J^{n-1}=\sqrt{|g|}:=\sqrt{\det[g_{ij}]}.
\end{equation}
 Since
 $\tau_t: S_x^{n-1}\to S_{\gamma_{\xi}(t)}^{n-1}$ is an
isometry, we have
$$
\langle d(\exp_x)_{t\xi}(t\eta_{a}),
d(\exp_x)_{t\xi}(t\eta_{b})\rangle_{g} =\langle
\mathbb{A}(t,\xi)(\eta_{a}), \mathbb{A}(t,\xi)(\eta_{b})\rangle_{g},
$$
and then $ \sqrt{|g|}=\det\mathbb{A}(t,\xi).$ So, by applying
(\ref{2.1}) and (\ref{2.2}), the volume $\mathrm{vol}(B(x,r))$ of a
geodesic ball $B(x,r)$, with radius $r$ and center $x$, on $M$ is
given by
\begin{eqnarray}\label{2.3}
\mathrm{Vol}(B(x,r))=\int\limits_{S_{x}^{n-1}}\int\limits_{0}^{\min\{r,d_{\xi}\}}\sqrt{|g|}dtd\sigma
=\int\limits_{S_{x}^{n-1}}\left(\int\limits_{0}^{\min\{r,d_{\xi}\}}\det(\mathbb{A}(t,\xi))dt\right)
d\sigma,
\end{eqnarray}
where $d\sigma$ denotes the $(n-1)$-dimensional volume element on
$\mathbb{S}^{n-1}\equiv S_{x}^{n-1}\subseteq{T_{x}M}$. As in Section
1, let $r(z)=d(x,z)$ be the intrinsic distance to the point
$x\in{M}$. Since for any $\xi\in{S}_{x}^{n-1}$ and $t_{0}>0$, we
have $\nabla{r}{(\gamma_{\xi}(t_{0}))}=\gamma'_{\xi}(t_{0})$ when
the point $\gamma_{\xi}(t_{0})=\exp_{x}(t_{0}\xi)$ is away from the
cut locus of $x$ (cf. \cite{a2}), then, by the definition of a
non-zero tangent vector ``\emph{radial}" to a prescribed point on a
manifold given in the first page of \cite{KK}, we know that for
$z\in{M}\backslash(Cut(x)\cup{x})$ the unit vector field
\begin{eqnarray*}
v_{z}:=\nabla{r}{(z)}
\end{eqnarray*}
is the radial unit tangent vector at $z$. We also need the following
fact about $r(z)$ (cf. Prop. 39 on p. 266 of \cite{pp}),
\begin{eqnarray*}
\partial_{r}\Delta{r}+\frac{(\Delta{r})^2}{n-1}\leq\partial_{r}\Delta{r}+|\mathrm{Hess}r|^{2}=-\mathrm{Ric}(\partial_{r},\partial_{r}),
\qquad {\rm{with}}~~\Delta{r}=\partial_{r}\ln( \sqrt{|g|}),
\end{eqnarray*}
  with
$\partial_{r}=\nabla{r}$ as a differentiable vector (cf. Prop. 7 on
p. 47 of \cite{pp}
 for the differentiation of $\partial_{r}$), where $\Delta$ is the Laplace operator on $M$ and $\mathrm{Hess}r$ is the Hessian of $r(z)$. Then,
together with (\ref{2.2}), we have
\begin{eqnarray}
&& J''+\frac{1}{(n-1)}\mathrm{Ric}\left(\gamma'_{\xi}(t),
\gamma'_{\xi}(t)\right)J\leq 0,  \nonumber\\
&&J(t,\xi)=t + O(t^2), \quad  J'(t,\xi)=1+O(t). \label{2.5}
\end{eqnarray}

\subsection{A volume comparison theorem in smooth metric measure spaces}

We also need the following volume comparison theorem proven by Wei
and Wylie (cf. \cite[Theorem 1.2]{ww}) which is the key point to
prove Theorem \ref{theorem1.5}.

\begin{theorem} (\cite{ww}) \label{theorem3.1}
Let $(M,g,e^{-f}dv_{g})$ be $n$-dimensional ($n\geqslant2$) complete
smooth metric measure space with $\mathrm{Ric}_{f}\geqslant(n-1)H$.
Fix $x_{0}\in{M}$. If $\partial_{t}f\geqslant-a$ along all minimal
geodesic segments from $x_{0}$ then for $R\geqslant{r}>0$ (assume
$R\leqslant\pi/2\sqrt{H}$ if $H>0$),
\begin{eqnarray*}
\frac{\mathrm{Vol}_{f}[B(x_{0},R)]}{\mathrm{Vol}_{f}[B(x_{0},r)]}\leqslant{e^{aR}}\frac{\mathrm{Vol}_{H}^{n}(R)}{\mathrm{Vol}_{H}^{n}(r)},
\end{eqnarray*}
where $\mathrm{Vol}_{H}^{n}(\cdot)$ is the volume of the geodesic
ball with the prescribed radius in the space $n$-form with constant
sectional curvature $H$, and, as before, $\mathrm{vol}_{f}(\cdot)$
denotes the weighted (or $f$-)volume of the given geodesic ball on
$M$. Moreover, equality in the above inequality holds if and only if
the radial sectional curvatures are equal to $H$ and
$\partial_{t}f\equiv-a$. In particular, if $\partial_{t}f\geqslant0$
and $\mathrm{Ric}\geqslant0$, then $M$ has $f$-volume growth of
degree at most $n$.
\end{theorem}

\subsection{Proof of Theorem \ref{theorem1.5}}

{\emph{Proof}}. For complete and \emph{non-compact} smooth metric
measure $n$-space $(M,g,e^{-f}dv_{g})$, if $\partial_{t}f\geqslant0$
(along all minimal geodesic segments from $x_{0}$) and
$\mathrm{Ric}_{f}\geqslant0$, then by Theorem \ref{theorem3.1} we
have
\begin{eqnarray*}
\frac{\mathrm{Vol}_{f}[B(x_{0},R)]}{\mathrm{Vol}_{f}[B(x_{0},r)]}\leqslant{e^{0\cdot{R}}}\cdot\frac{V_{0}(R)}{V_{0}(r)}=\left(\frac{R}{r}\right)^{n},
\end{eqnarray*}
with, as before, $V_{0}(\cdot)$ denotes the volume of the ball with
the prescribed radius in $\mathbb{R}^{n}$. Clearly, here the volume
doubling condition (\ref{a2}) is satisfied with $C_0=1$.

For $(M,g,e^{-f}dv_{g})$, in order to apply Theorem \ref{theorem1.2}
to prove Theorem \ref{theorem1.5}, we need to normalize the original
measure $e^{-f}dv_{g}$ such that the volume condition (\ref{a3}) can
be satisfied. In fact, we need to choose the positive measure $\mu$
 to be $\mu=e^{f(x_0)-f}dv_{g}$. Then by applying (\ref{wv}), (\ref{2.2}), (\ref{2.3}) and
(\ref{2.5}), we can get
\begin{eqnarray*}
\lim_{r\rightarrow 0}\frac{\mu (B(x_0,r))}{\mu_{E}
(\mathbb{B}_{n}(r))}=\lim_{r\rightarrow 0}\frac{e^{f(x_0)}\cdot
\mathrm{Vol}_{f}[B(x_{0},r)]}{V_{0}(r)}&=& \lim_{r\rightarrow
0}\frac{\int\limits_{\mathbb{S}^{n-1}}\left(\int\limits_{0}^{\min\{r,d_{\xi}\}}J^{n-1}(t,\xi)\cdot{e}^{-f}dt\right)d\sigma}
{e^{-f(x_0)}\int\limits_{\mathbb{S}^{n-1}}\int\limits_{0}^{r}t^{n-1}dtd\sigma}\\
&=&\frac{J'(0,\xi)\cdot{e}^{-f(x_{0})}}{e^{-f(x_0)}}=1
\end{eqnarray*}
by applying L'H\^{o}pital's rule $n$-times, which implies (\ref{a3})
is satisfied. Therefore, if in addition the Gagliardo-Nirenberg type
inequality (\ref{a5}) is satisfied, then by applying Theorem
\ref{theorem1.2}, we can get (\ref{a6}) directly. This completes the
proof of Theorem \ref{theorem1.5}.   \hfill $\Box$

\section*{Acknowledgments}
\renewcommand{\thesection}{\arabic{section}}
\renewcommand{\theequation}{\thesection.\arabic{equation}}
\setcounter{equation}{0} \setcounter{maintheorem}{0}

F. Du and J. Mao were partially supported by the NSF of China (Grant
No. 11401131) and CNPq, Brazil.

\end{document}